\documentclass{article}

\def\qed{\ifmmode\square\else\nolinebreak\hfill$\diamondsuit$\fi\par\vskip12pt}
\newcommand{\be}{\begin{equation}}
\newcommand{\ee}{\end{equation}}

\newcommand{\beq}{\begin{eqnarray}}
\newcommand{\eeq}{\end{eqnarray}}
\newcommand{\nbeq}{\begin{eqnarray*}}
\newcommand{\neeq}{\end{eqnarray*}}

\newcommand{\AmS}{{\protect\the\textfont2
  A\kern-.1667em\lower.5ex\hbox{M}\kern-.125emS}}

\def\be{\begin{equation}}
\def\ee{\end{equation}}

\begin{document}

\title
{\bf Characterizations of Student's t-distribution via regressions of order statistics}%

\author{George P. Yanev$^a$\thanks{Correspponding author. Email: yanevgp@ytpa.edu} \ and  M. Ahsanullah$^b$ \\
$^a$The University of Texas - Pan American, Edinburg, Texas, USA \\
$^b$Rider University, Lawrenceville, New Jersey, USA }%

\maketitle

\begin{abstract}
\noindent Utilizing regression properties of order statistics, we characterize
a family of distributions introduced by
Akhundov et al. \cite{ABN04}, that includes the t-distribution with two degrees of freedom as one of its members.
Then we extend this characterization result to t-distribution with more than two degrees of
freedom.


\end{abstract}%

%\begin{center}

%\end{center}

\section{Discussion of the Results}

The Student's $t$-distribution is widely used in statistical inferences when
the population standard deviation is unknown and is substituted by its
estimate from the sample. Recently Student's $t$-distribution was also
considered in financial modeling by Ferguson and Platen \cite{FP06} and as a
pedagogical tool by Jones \cite{J02}. The probability density function (pdf) of
the $t$-distribution with $\nu$ degrees of freedom ($t_{\nu}$-distribution) is
given for $-\infty<x<\infty$ and $\nu=1,2\ldots$ by
\begin{equation}
\label{pdf}f_{\nu}(x)=c_{\nu}\left(  1+\frac{x^{2}}{\nu}\right)  ^{-(\nu+1)/2}
\qquad\mbox{where}\quad c_{\nu}=\frac{\Gamma\left(  \frac{\nu+1}{2}\right)
}{\Gamma\left(  \frac{\nu}{ 2}\right)  \sqrt{\pi\nu} }%
\end{equation}
and $\Gamma(x)$ is the gamma function.

The vast majority of characterization results for univariate continuous
distributions based on ordered random variables is concentrated to exponential
and uniform families. It was not until recently, last 7-8 years, when some characterizations were obtained for 
$t_\nu$-distribution with $\nu=2$ and $\nu=3$. In this note we communicate generalizations of these recent
results for $t_{\nu}$-distribution when $\nu\ge 2$. Let $X, X_{1}, X_{2}, \ldots,
X_{n}$ for $n\ge3$ be independent random variables with common  cumulative
distribution function (cdf) $F(x)$. Assume that $F(x)$ is absolute continuous
with respect to the Lebesgue measure. Let $X_{1:3}\le X_{2:3}\le\ldots\le
X_{n:n}$ be the corresponding order statistics. Nevzorov et al. \cite{NBA03} (see
also Nevzorov \cite{N02}) and Akhundov and Nevzorov \cite{AN10} prove characterizations
for the $t_{\nu}$-distribution when $\nu= 2$ and $\nu=3$, respectively,
assuming, in addition,  $n=3$. Here we extend these results to the general case of any
$\nu\ge 2$ and any $n\ge3$.

Let $Q(x)$ be the quantile function of a random variable with cdf $F(x)$,
i.e., $F(Q(x))=x$ for $0<x<1$. Akhundov et al. \cite{ABN04} prove that for
$0<\lambda<1$ the relation
\begin{equation}
\label{ABN04}E[\lambda X_{1:3}+(1-\lambda)X_{3:3}\ |\ X_{2:3}=x]=x
\end{equation}
characterizes a family of probability distributions with quantile function
\begin{equation}
\label{family}Q_{\lambda}(x)=\frac{c(x-\lambda)}{\lambda(1-\lambda
)(1-x)^{\lambda}x^{1-\lambda}}+d, \qquad0<x<1,
\end{equation}
where $0<c<\infty$ and $-\infty<d<\infty$. Let us call this family of
distributions - $Q$-family.

{\bf\  Theorem 1 ($Q$-family)}\ {\it Assume that $E|X|<\infty$ and $n\ge3$
is a positive integer. The random variable $X$ belongs to the $Q$-family if
and only if for some $2\le k\le n-1$ and some $0<\lambda<1$
\beq
\label{eqn_one_t2_b}\lefteqn{\lambda E\left[  \frac{1}{k-1}\sum_{i=1}%
^{k-1}(X_{k:n}-X_{i:n})\ |\ X_{k:n}=x\right]  }  &  & \nonumber\\
&  & = (1-\lambda) E\left[  \frac{1}{n-k}\sum_{j=k+1}^{n}(X_{j:n}%
-X_{k:n})\ |\ X_{k:n}=x\right]  .
\eeq
} Note that, (\ref{eqn_one_t2_b}) can be written as
\begin{equation}
\label{eqn_one_t23}\lambda E\left[ \frac{1}{k-1}\sum_{j=1}^{k-1}%
X_{j:n}\ |\ X_{k:n} =x\right] +(1-\lambda) E\left[ \frac{1}{n-k}\sum
_{j=k+1}^{n}X_{j:n}\ |\ X_{k:n}=x\right] =x.
\end{equation}
Clearly for $n=3$ and $k=2$, (\ref{eqn_one_t23})  reduces to (\ref{ABN04}). It
is also worth mentioning here that, as Balakrishnan and Akhundov \cite{BA03}
report, the $Q$-family, for different values of $\lambda$, approximates well a
number of common distributions including Tukey lambda, Cauchy, and Gumbel (for maxima).

Notice that $t_{2}$-distribution belongs to the $Q$-family, having quantile
function (e.g., Jones \cite{J02})
\[
Q_{1/2}(x)=\frac{2^{1/2}(x-1/2)}{x^{1/2}(1-x)^{1/2}}, \qquad0<x<1.
\]
Nevzorov et al. (2003) prove that if $E|X|<\infty$ then $X$ follows $t_{2}%
$-distribution if and only if
\begin{equation}
\label{N03}E[X_{2:3}-X_{1:3}|X_{2:3}=x]= E[X_{3:3}-X_{2:3}|X_{2:3}=x].
\end{equation}
This also follows directly from (\ref{ABN04}) with $\lambda=1/2$. Recall that
the cdf of $t_{2}$-distribution (see Jones \cite{J02}) is
\[
%be \label{F2}
F_{2}(x)=\frac{1}{2}\left(  1+\frac{x}{\sqrt{1+x^{2}}}\right)  .
\]
%ee
Setting $\lambda=1/2$ in (\ref{eqn_one_t2_b}), we obtain the following
corollary of Theorem 1.

{\bf\  Corollary ($t_{2}$-distribution)}\ {\it\ Assume that $E|X|<\infty$
and $n\ge3$ is a positive integer. Then
\begin{equation}
\label{cdf_t2}F(x)=F_{2}\left(  \frac{x-\mu}{\sigma}\right)  \qquad\mbox{for}
\quad -\infty<\mu<\infty, \quad\sigma>0,
\end{equation}
if and only if for some $2\le k\le n-1$
\be
\label{eqn_one_t2}E\left[  \frac{1}{k-1}\sum_{i=1}^{k-1}%
(X_{k:n}-X_{i:n})\ |\ X_{k:n}=x\right] 
=E\left[  \frac{1}{n-k}\sum_{j=k+1}^{n}(X_{j:n}-X_{k:n})\ |\ X_{k:n}%
=x\right]  .
\ee}

\noindent Relation (\ref{eqn_one_t2}) can be interpreted as follows. Given the value
of $X_{k:n}$, the average deviation from $X_{k:n}$ to the observations less than it
equals the average deviation from the observations greater than $X_{k:n}$
to  it.

{\bf\ Remarks}\ (i) Notice that (\ref{eqn_one_t2})
reduces to (\ref{N03}) when $n=3$ and $k=2$. (ii) Let us set $n=2r+1$ and
$k=r+1$ for an integer $r\ge 1$. Let $M_{2r+1}=X_{r+1:2r+1}$ be the median of the sample $X_1, X_2,\ldots , X_{2r+1}$. Then (\ref{eqn_one_t2}) implies
\[
E\left[ \sum_{i=1}^{r}%
(M_{2r+1}-X_{i:2r+1})|M_{2r+1}=x\right]\!
=E\! \left[ \sum_{j=r+2}^{2r+1}(X_{j:2r+1}-M_{2r+1})|M_{2r+1}
=x\right]  .
\]
If, in addition, $\bar{X}_{2r+1}=\sum_{i=1}^{2r+1}X_i/(2r+1)$ is the sample mean, then (\ref{eqn_one_t2}) reduces to Nevzorov et al. \cite{NBA03} $t_{2}%
$-distribution characterization relation
\[
E\left[ \bar{X}_{2r+1}\ |\ M_{2r+1}=x\right]  =x.
\]

Let us now turn to the case of $t_{\nu}$-distribution with $\nu\ge3$. Akhundov
and Nevzorov \cite{AN10} extend (\ref{N03}) to a characterization of $t_{3}%
$-distribution as follows. If $EX^{2}<\infty$ then $X$ follows $t_{3}%
$-distribution if and only if
\begin{equation}
\label{AN10}E[(X_{2:3}-X_{1:3})^{2}|X_{2:3}=x]=E[(X_{3:3}-X_{2:3})^{2}%
|X_{2:3}=x].
\end{equation}
We generalize this in two directions: (i) characterizing $t_{\nu}%
$-distribution with $\nu\ge3$ and (ii) considering a sample of size $n\ge3$.
The following result holds.

{\bf\ Theorem 2 ($t_{\nu}$-distribution)}\ {\it\ Assume $EX^{2}<\infty$.
Let $n\geq3$ and $\nu\geq3$ be positive integers. Then
\begin{equation}
F(x)=F_{\nu}\left(  \frac{x-\mu}{\sigma}\right)  \qquad\mbox{for}\quad
-\infty<\mu<\infty,\quad\sigma>0,\label{F_nu}%
\end{equation}
where $F_{\nu}(x)$ is the $t_{\nu}$-distribution cdf if and only if for some
$2\leq k\leq n-1$
\beq
\lefteqn{E\left[  \frac{1}{k-1}\sum_{i=1}^{k-1}\left(  \frac{\nu-1}{2}%
X_{k:n}-(\nu-2)X_{i:n}\right)  ^{2}\ |\ X_{k:n}=x\right]  }%
\nonumber\label{new_eqn10}\\
&  = & E\left[  \frac{1}{n-k}\sum_{j=k+1}^{n}\left(  (\nu-2)X_{j:n}-\frac{\nu
-1}{2}X_{k:n}\right)  ^{2}\ |\ X_{k:n}=x\right]  .
\eeq}
{\bf\ Remarks}\ (i) Notice that if $n=3$, $k=2$, and $\nu=3$, then (\ref{new_eqn10}) reduces to
(\ref{AN10}). (ii) Let us set $\nu=3$, $n=2r+1$ and
$k=r+1$ for an integer $r\ge 1$. If, as before,  $M_{2r+1}=X_{r+1:2r+1}$ is the median of the sample $X_1, X_2,\ldots , X_{2r+1}$, then (\ref{new_eqn10}) implies the following equality between the sum of squares of the deviations from the sample median
\[
E\! \left[ \sum_{i=1}^{r}\!\left(
M_{2r+1}-X_{i:2r+1}\right)^{2}|M_{2r+1}=x\right] \!
=\!E\! \left[ \sum_{j=r+2}^{2r+1}\! \left( X_{j:2r+1}-M_{2r+1}\right)^{2}|M_{2r+1}=x\right].
\]

\section{ Proofs}

To prove our results we need the following two lemmas.

{\bf\ Lemma 1}\ (Balakrishnan and Akhundov \cite{BA03}) The cdf $F(x)$ of a
random variable $X$ with quintile function (\ref{family}) is the only
continuous cdf solution of the equation
\begin{equation}
\label{lemma1}[F(x)]^{2-\lambda}[1-F(x)]^{1+\lambda}=cF^{\prime}(x), \qquad
c>0.
\end{equation}

{\bf\ Lemma 2}\ Let $r\ge1$ and $n\ge2$ be integers. Then
\beq
\label{lemma 2}\frac{1}{k-1}\sum_{i=1}^{k-1}E\left[  X_{i:n}^{r}%
\ |\ X_{k:n}=x\right]   &  = & \frac{1}{F(x)}\int_{-\infty}^{x} t^{r} dF(t),
\qquad 2\le k\le n;\\
\frac{1}{n-k}\sum_{j=k+1}^{n}E\left[  X_{j:n}^{r}\ |\ X_{k:n}=x\right]   &  = & 
\frac{1}{1-F(x)}\int_{x}^{\infty}t^{r} dF(t), \qquad 1\le k\le n-1.\nonumber
\eeq

{\bf\ Proof.}\ Using the standard formulas for the conditional density of
$X_{j:n}$ given $X_{k:n}=x$ $(j<k)$ (e.g., Ahsanullah and Nevzorov \cite{AN01},
Theorem 1.1.1), we obtain for $r\ge 1$
\nbeq
\lefteqn{\frac{1}{k-1}\sum_{j=1}^{k-1}E[X_{j:n}^{r}\ |\ X_{k:n}=x]}\\
&  = & \frac{1}{(k-1)}\frac{(k-1)}{[F(x)]^{k-1}}\sum_{j=1}^{k-1}{k-2
\choose j-1}\int_{-\infty}^{x}[F(t)]^{j-1}[F(x)-F(t)]^{k-1-j}t^{r}dF(t)\\
&  = & \frac{1}{[F(x)]^{k-1}}\sum_{i=0}^{k-2}{k-2 \choose i}\int_{-\infty}%
^{x}[F(t)]^{i}[F(x)-F(t)]^{k-2-i}t^{r}dF(t)\\
&  = & \frac{1}{F(x)}\int_{-\infty}^{x}t^{r}dF(t).
\neeq
This verifies (\ref{lemma 2}). The second relation in the lemma's statement
can be proved similarly.

\subsection{Proof of Theorem 1}

First, we show that equation (\ref{eqn_one_t2_b}) implies (\ref{family}).
Applying Lemma 2, for the left-hand side of (\ref{eqn_one_t23}), we obtain
\be
\label{eqn_one_t24}\frac{\lambda}{k-1}\sum_{j=1}^{k-1}E[X_{j:n}
\ |\ X_{k:n}=x]+\frac{1-\lambda}{n-k}\sum_{j=k+1}^{n}E[X_{j:n}\ |\ X_{k:n}%
=x] = \frac{\lambda}{F(x)}\int_{-\infty}^{x}tdF(t) + \frac{1-\lambda}%
{1-F(x)}\int_{x}^{\infty}tdF(t).
\ee

Further, since $E|X|<\infty$, we have
\begin{equation}
\label{limits}\lim_{x\to-\infty}xF(x)=0\qquad\mbox{and}\qquad\lim_{x\to\infty
}x(1-F(x))=0.
\end{equation}
Therefore, integrating by parts, we obtain
\be
\label{int_parts}\frac{\lambda}{F(x)}\int_{-\infty}^{x}%
tdF(t)+\frac{1-\lambda}{1-F(x)}\int_{x}^{\infty}tdF(t)
= x-\frac{\lambda}{F(x)}\int_{-\infty}^{x} F(t)dt+\frac{1-\lambda
}{1-F(x)}\int_{x}^{\infty}(1-F(t))dt.
\ee 
Thus, from (\ref{eqn_one_t24}) and (\ref{int_parts}) it follows that
(\ref{eqn_one_t2_b}) is equivalent to
\[
\label{main}\lambda(1-F(x))\int_{-\infty}^{x} F(t)dt=(1-\lambda) F(x)\int
_{x}^{\infty}(1-F(t))dt .
\]
The last equation can be written as
\[
-\frac{\lambda}{1-\lambda}\int_{-\infty}^{x}\!\!F(t)dt\ \frac{d}{dx} \left[
\int_{x}^{\infty}\!\! (1-F(t))dt\right]  =\int_{x}^{\infty}\!\!
(1-F(t))dt\ \frac{d}{dx} \left[  \int_{-\infty}^{x} \! \! F(t)dt\right]  ,
\]
which leads to
\[
\int_{-\infty}^{x}F(t)dt=c \left(  \int_{x}^{\infty}(1-F(t))dt\right)
^{-\lambda/(1-\lambda)}\qquad c>0.
\]
Differentiating both sides with respect to $x$ we obtain
\[
\int_{x}^{\infty}(1-F(t))dt=c_{1}\left(  \frac{1}{F(x)}-1\right)  ^{1-\lambda
}, \qquad c_{1}>0.
\]
Differentiating one more time, we have
\begin{equation}
\label{lemma1a}[F(x)]^{2-\lambda}[(1-F(x))]^{1+\lambda}=c_{2}F^{\prime}(x),
\qquad c_{2}>0,
\end{equation}
which is (\ref{lemma1}). Referring to Lemma 1 we see that (\ref{eqn_one_t2_b})
implies (\ref{family}).

To complete the proof of the theorem, it remains to verify that $F(x)$ with
quantile function (\ref{family}) satisfies (\ref{eqn_one_t2_b}).
Differentiating (\ref{family}) with respect to $x$ we obtain
\[
Q^{\prime}_{\lambda}(x)=c(1-x)^{-(1+\lambda)}x^{-(2-\lambda)}\qquad c>0.
\]
On the other hand, since $F(Q_{\lambda}(x))=x$, we have $Q^{\prime}_{\lambda
}(x)=[F^{\prime}(Q_{\lambda}(x))]^{-1}$. Therefore,
\[
(1-x)^{1+\lambda}x^{2-\lambda}=cF^{\prime}(Q_{\lambda}(x)),
\]
which is equivalent to (\ref{lemma1a}) and thus, to (\ref{eqn_one_t2_b}). This
completes the proof.

\subsection{Proof of Theorem 2}

Notice that (\ref{new_eqn10}) can be written as
\nbeq
%\label{new_eqn2}
\lefteqn{(\nu-1)x\left[\frac{1}{n-k}\sum_{j=k+1}^{n}E[X_{j:n}\ |\ X_{k:n}=x]- \frac{1}{k-1}\sum_{j=1}^{k-1}E[X_{j:n}
\ |\ X_{k:n}=x]\right]}  \\
&  & =(\nu-2)\left[  \frac{1}{n-k}\sum_{j=k+1}^{n}E[X_{j:n}^{2}\ |\ X_{k:n}=x]- \frac{1}{k-1}\sum_{j=1}^{k-1}E[X_{j:n}^{2}\ |\ X_{k:n}=x]\right].
\neeq

Referring to Lemma 2 with $r=1$ and $r=2$, we see that this is
%(\ref{new_eqn2})
equivalent to
\beq
\label{intermediate_eqn}\lefteqn{(\nu-1)x\left[  \frac{1}{1-F(x)}\int
_{x}^{\infty}tdF(t)- \frac{1}{F(x)}\int_{-\infty}^{x} tdF(t)\right]  }  &  &
\nonumber\\
&  & =(\nu-2)\left[  \frac{1} {1-F(x)}\int_{x}^{\infty}t^{2}dF(t)-\frac
{1}{F(x)}\int_{-\infty}^{x} t^{2}dF(t)\right]  .
\eeq
Let us assume that $EX=0$ and $EX^{2}=1$. Hence
\[
\int_{x}^{\infty}tdF(t)=-\int_{-\infty}^{x} tdF(t)\quad\mbox{and}\quad\int
_{x}^{\infty}t^{2}dF(t)=1-\int_{-\infty}^{x} t^{2}dF(t)
\]
and thus (\ref{intermediate_eqn}) is equivalent to
\[
-(\nu-1)x\left(  \frac{1}{1-F(x)}+\frac{1}{F(x)}\right)  \int
_{-\infty}^{x} tdF(t)
= \frac{\nu-2}{1-F(x)}-(\nu-2)\left(  \frac{1}{1-F(x)}+\frac{1}{F(x)}\right)  \int_{-\infty
}^{x} t^{2}dF(t)
\]
Multiplying the above equation by $F(x)[1-F(x)]$, we find
\begin{equation}
\label{star}-(\nu-1)x\int_{-\infty}^{x} tdF(t)=(\nu-2)\left[  F(x) -
\int_{-\infty}^{x} t^{2}dF(t)\right]  .
\end{equation}
Differentiating both sides with respect to $x$, we obtain
\[
-(\nu-1)\int_{-\infty}^{x} tdF(t)=f(x)(x^{2}+\nu-2).
\]
Since the left-hand side of the above equation is differentiable, we have
that $f^{\prime}(x)$ exists. Differentiating both sides with respect to $x$,
we find
\[
\frac{f^{\prime}(x)}{f(x)}=-\frac{\nu+1}{\nu-2}\ \frac{x}{1+\frac
{\displaystyle x^{2}}{\displaystyle \nu-2}}.
\]
Integrating both sides and making use of the fact that $f(x)$ is a pdf we
obtain
\begin{equation}
\label{Zpdf}f(x)=c\left(  1+\frac{x^{2}}{\nu-2}\right)  ^{-(\nu+1)/2}%
\quad\mbox{where}\quad c=\frac{\Gamma\left(  \frac{\nu+1}{2}\right)  }%
{\Gamma\left(  \frac{\nu}{ 2}\right)  \sqrt{(\nu-2)\pi} }.
\end{equation}
It is not difficult to see that if a random variable $Z$ has the pdf
(\ref{Zpdf}), then
\[
X=Z\sqrt{\frac{\nu}{\nu-2}}
\]
follows $t_{\nu}$-distribution, i.e., its pdf is given by (\ref{pdf}). Thus,
we have proved that (\ref{new_eqn10}) implies (\ref{F_nu}) when $\mu=0$ and
$\sigma^{2}=1$. The result now follows in the general case by considering the
linear transformation $Y=\sigma X+\mu$.

To complete the proof, we need to verify that (\ref{new_eqn10}) holds when $X$
has a cdf given by (\ref{F_nu}). If $X$ has pdf (\ref{pdf}) (i.e., cdf
(\ref{F_nu})), then we define
\[
Z=X\sqrt{\frac{\nu-2}{\nu}},
\]
which has pdf (\ref{Zpdf}). Now, it is not difficult to verify that
(\ref{Zpdf}) satisfies (\ref{star}), which in turn is equivalent to
(\ref{new_eqn10}). The proof is complete.


\begin{thebibliography}{9}                                                                                                %
%{xbib}


\bibitem {AN01}Ahsanullah, M. and Nevzorov, V.B. (2001). Ordered Random
Variables, NOVA Sci. Publ., Huntington, NY.

\bibitem {ABN04}Akhundov, I.S., Balakrishnan, N., Nevzorov, V.B. (2004). New
characterizations by properties of midrange and related statistics.
Communications in Statistics: Theory and Methods, 33(12), 3133-3143.

\bibitem {AN10}Akhundov, I. and Nevzorov, V.B. (2010). A simple
characterization of Student's $t_{3}$ distribution. Stats. Probab. Lett.,
80(5-6), 293-295.

\bibitem {BA03}Balakrishnan, N. and Akhundov, I. (2003). A characterization by
linearity of the regression function based on order statistics. Stats. Probab.
Lett., 63(4), 435-440.

\bibitem {FP06}Ferguson, K. and Platen, E. (2006). On the distributional
characterization of daily log-returns of a world stock index. Appl.
Math. Finance, 13(1), 19-38.

\bibitem {J02}Jones, M.C. (2002). Student's simplest distribution. The Statististician,
51(1), 41-49.

\bibitem {N02}Nevzorov, V.B. (2002). On a property of Student's distribution
with two degrees of freedom. Zapiski Nauchnykh Seminarov POMI, Vol. 294,
148-157 (in Russian). (English Translation: J. Math. Sci., 2005, 127(1), 1757-1762.)

\bibitem {NBA03}Nevzorov, V.B., Balakrishnan, N., and Ahsanullah, M. (2003).
Simple characterizations of Student's $t_{2}$-distribution. The Statistician,
52(3), 395-400.
\end{thebibliography}
\end{document}